\documentclass[11pt]{article}
 \usepackage{amsmath,amssymb}
 \usepackage{epsfig}
 \usepackage{graphicx}
 \usepackage{pict2e}
\usepackage{verbatim}
\usepackage{multirow} 
\usepackage{bm}
 \newtheorem{Lemma}{Lemma}
 \newtheorem{Proposition}[Lemma]{Proposition}
 
 \newtheorem{Theorem}[Lemma]{Theorem}

 \newcommand{\FF}{\mbox{${\mathcal F}$}} 
  \newcommand{\EE}{\mbox{${\mathcal E}$}}

 \newcommand{\sfrac}[2]{{\textstyle\frac{#1}{#2}}}

  \newcommand{\Ints}{{\mathbb{Z}}}
  \newcommand{\bw}{{\mathbf w}}
    \newcommand{\bV}{{\mathbf V}}
  \newcommand{\bE}{{\mathbf E}}
    
 \newcommand{\bt}{{\mathbf t}}    
 
    \newcommand{\eps}{\varepsilon}

\def\Ex{\mathbb{E}}
\def\Pr{\mathbb{P}}
\def\ind{{\rm 1\hspace{-0.90ex}1}}

\newcommand{\var}{\mathrm{var}}

\newcommand{\Ttria}{T^{tria}}
\newcommand{\Tspan}{T^{span}}

\begin{document}

\title{Weak Concentration for First Passage Percolation Times on Graphs and General Increasing Set-valued Processes}
 \author{David J. Aldous
 \thanks{Department of Statistics,
 367 Evans Hall \#\  3860,
 U.C. Berkeley CA 94720;  aldous@stat.berkeley.edu;
  www.stat.berkeley.edu/users/aldous.  Aldous's research supported by
 N.S.F Grant DMS-106998. }}

 \maketitle
 
  \begin{abstract}
A simple lemma bounds $\mathrm{s.d.}(T)/\Ex T$ for hitting times $T$ in Markov chains with a certain strong monotonicity property. 
We show how this lemma may be applied to several increasing set-valued processes.
Our main result concerns  a model of first passage percolation on a finite graph, where the traversal times of edges are independent Exponentials with arbitrary rates.  
Consider the percolation time $X$ between two arbitrary vertices. 
We prove that $\mathrm{s.d.}(X)/\Ex X$ is small if and only if $\Xi/\Ex X$ is small, where $\Xi$ is the maximal edge-traversal time
in the percolation path attaining $X$.  
  \end{abstract}
  
 \vspace{0.1in}

 \vspace{0.1in}
 {\em Key words and phrases;} Concentration inequality, monotone, Markov chain, coupling, first passage percolation.

\section{Introduction}
\label{sec:intro}
For random variables $T$ arising (loosely speaking) as optimal costs in some ``random environment" model,
one might not be able to estimate $\Ex T$ explicitly, for one of two reasons:
it may involve a difficult optimization problem (exemplified by the Euclidean TSP over random points \cite{steele}), 
or the model may involve many parameters (the case more relevant to this paper).
In such cases the well-known {\em method of bounded differences}  \cite{mcdiarmid}
often enables us to bound $T - \Ex T$ explicitly, 
and this general topic of {\em concentration inequalities} has been developed in many directions 
over the last generation.
What we will call (recalling the {\em weak} law of large numbers) 
a {\em weak  concentration inequality} 
is just a result showing  $\mathrm{s.d.}(T)/\Ex T$ is small.
The starting point for this paper is a simple technique (section \ref{sec:monot}) 
for proving  weak  concentration  for hitting times $T$ in Markov chains with a certain strong monotonicity property. 
At first sight one might doubt that this technique could be more than very narrowly applicable, but its quick use in 
diverse other contexts (sections \ref{sec:growth} - \ref{sec:coverage}; see also section \ref {sec:remarks}(d)) should partly dispel such doubts. 
Our main result concerns a more sophisticated analysis of 
a model of first passage percolation on a finite graph, where the traversal times of edges are independent Exponentials with arbitrary rates.  
We will  describe that result and its background  in section \ref{sec:introperc}, and prove it in the subsequent sections.

.

\subsection{A monotonicity condition}
\label{sec:monot}
The setting is a continuous-time Markov chain $(Z_t)$ on a finite state space $\Sigma$, 
where we study the hitting time 
\begin{equation}
 T := \inf \{t: \ Z_t \in \Sigma_0\} 
 \label{T-def}
 \end{equation}
for a fixed subset $\Sigma_0 \subset \Sigma$.  
Assume
\begin{equation}
\label{h-def}
 h(S) := \Ex_S T < \infty \mbox{ for each } S \in \Sigma 
 \end{equation}
which holds under the natural ``reachability" condition.
Assume also a rather strong ``monotonicity" condition:
\begin{equation}
h(S^\prime) \le h(S) \mbox{ whenever $S \to S^\prime$ is a possible transition}. 
\label{def-monotone}
\end{equation}
In a typical example, the state space $\Sigma$ will be the space of all subsets $S$ of a given finite set $\bV$, and possible transitions 
will be of the form $S \to S \cup \{v\}$.
A special case of our simple lemma, which suffices for the quick applications in sections \ref{sec:growth} - \ref{sec:coverage}, is
\begin{Lemma}
\label{L-simple-1}
Under condition (\ref{def-monotone}), for any initial state,
\[ \frac{\var  \ T}{\Ex  T} \le  \max \{ h(S) - h(S^\prime)    : \ S \to S^\prime \mbox{ a possible transition} \} .
\]
\end{Lemma}
In the ``quick applications" below, the right side will be bounded using natural couplings of the processes started from $S$ and from $S^\prime$.
For more complicated applications, Lemma \ref{L-simple-2} allows  the possibility of  occasional transitions for which 
$ h(S) - h(S^\prime) $ is not so small.
\begin{Lemma}
\label{L-simple-2}
In the setting of Lemma \ref{L-simple-1}, for arbitrary $\delta, \eps > 0$,
\[
 \frac{\var \ T}{(\Ex T)^2} \le \delta  + \eps
+\frac{  \Ex \int_0^T \ind_{\{q_\delta (Z_u) \ge \eps\}} du }{\Ex T} . 
\]
where $ q_\delta (S)$ is defined at (\ref{def-q}) below.
\end{Lemma}
The proofs below involve only very standard martingale analysis. 
We are not claiming any conceptual novelty in these results, but instead emphasize their applications later. 
A variant form of Lemma \ref{L-simple-2} is derived in \cite{me-incipient} Corollary 3.

{\bf Proof of Lemma \ref{L-simple-1}.}  
Expectation is relative to some fixed initial state $S_0$.
 Note that
 $T  = \inf \{t: h(Z_t) = 0\}$, that
 $t \to h(Z_t)$ is decreasing and so for any reachable state $S$ we have
 $h(S) \le h(S_0) = \Ex T$,  facts we use frequently without comment.  
 Consider the martingale 
 \begin{equation}
  M_t := \Ex (T \vert Z_t) = h(Z_{t \wedge T}) + {t \wedge T}. 
  \label{MgT}
  \end{equation}
 The Doob-Meyer decomposition of $M_t^2$ into a martingale $Q_t$ and a predictable process is clearly
 \[ M_t^2 - M_0^2 = Q_t + \int_0^t a(Z_u) \ du \]
 where
 \begin{equation}
  a(S) := \sum_{S^\prime} q(S,S^\prime) \ (h(S) - h(S^\prime) )^2 
 \label{def-a}
 \end{equation}
where $q(S,S^\prime) $ are the transition rates.  
 Taking expectation at $t=\infty$ gives 
 \begin{equation}
  \var \  T = \Ex  \int_0^T a(Z_u) \ du . 
  \label{Tv}
  \end{equation}
 On the other hand the martingale property (\ref{MgT}) for $\Ex (T\vert Z_t)$ corresponds to the identity
 \begin{equation}
  b(S) :=  \sum_{S^\prime} q(S,S^\prime) (h(S) - h(S^\prime) ) =  1  \mbox{ while } S \cap \Sigma_0 = \emptyset 
  \label{def-b}
  \end{equation}
 and therefore 
 \begin{equation}
 \Ex T =  \Ex \int_0^T b(Z_u) \ du . 
 \label{TE}
 \end{equation}
 Setting 
 \[ \kappa = \max \{ h(S) - h(S^\prime)    : \ S \to S^\prime \mbox{ a possible transition} \}  \]
 we clearly have $a(S) \le \kappa b(S)$; 
 note this is where we use the monotonicity hypothesis (\ref{def-monotone}).
The result now follows from (\ref{Tv},\ref{TE}).

 {\bf Proof of Lemma \ref{L-simple-2}.}  
 We continue with the notation above.
Fix $\delta > 0$ and write, for a possible transition $S \to S^\prime$,  
\[ h(S) - h(S^\prime) \le \delta\ \Ex T 
+ \ind_{\{h(S) - h(S^\prime) > \delta \Ex  T\} }
 \ \Ex T .
\]
Using this to bound one term of the product 
$(h(S) - h(S^\prime)^2$ in the definition (\ref{def-a}) of $a(S)$, 
and comparing with  the definition (\ref{def-b}) of $b(S)$, we obtain
\begin{eqnarray*}
a(S) & \le & b(S) \ \cdot \ \delta \Ex T \\
&& + \sum_{S^\prime : \ h(S) - h(S^\prime) >  \delta \Ex T } 
\ q(S,S^\prime) (h(S) - h(S^\prime) )  \ \cdot \ \Ex T
\end{eqnarray*}
While  $S \cap \Sigma_0 = \emptyset$  we have $b(S) = 1$ and so
\[ \frac{a(S)}{  \Ex T }
\le \delta +  q_\delta (S) \]
where
\begin{equation}
 q_\delta (S) :=  \sum_{S^\prime : \ h(S) - h(S^\prime) >  \delta \Ex T } 
\ q(S,S^\prime) (h(S) - h(S^\prime) ) . 
\label{def-q}
\end{equation}
Using (\ref{Tv})
\[ \frac{\var \ T}{\Ex T} \le \delta \Ex T
+ \Ex \int_0^T q_\delta (Z_u) du . \]
Because $q_\delta (S) \le b(S) \equiv 1$ we can fix $\eps > 0$ and write
$q_\delta (Z_u) \le \eps + \ind_{\{q_\delta (Z_u) \ge \eps\}}$,
and the result follows.

\paragraph{Discrete-time chains.}
Given a discrete-time Markov chain with transition probabilities $p(S,S^\prime)$ there is a corresponding 
continuous-time Markov chain with transition rates  $q(S,S^\prime) = p(S,S^\prime)$.
The relation between the hitting times $ T_{\mathrm{disc}}$ and $T_{\mathrm{cont}}$ for these two chains is (using Wald's identity)
\begin{equation}
\Ex T_{\mathrm{cont}} = \Ex T_{\mathrm{disc}}; \quad 
\var \ T_{\mathrm{cont}} = \var \  T_{\mathrm{disc}} + \Ex T_{\mathrm{disc}} .
\label{continuization}
\end{equation}
Via this {\em continuization} device our results may be applied to discrete-time chains.

\subsection{A general Markovian growth process}
\label{sec:growth}
For a first ``obvious" use of Lemma \ref {L-simple-1},
we consider a general growth process $(Z_t)$ on the lattice $\Ints^2$.
The states are finite vertex-sets $S$, the possible transitions are 
$S \to S \cup \{v\}$ where $v$ is a vertex adjacent to $S$.  
For each such transition, we assume the transition rates are bounded above and below:
\begin{equation}
0 < c_* \le q(S, S \cup \{v\}) \le c^* < \infty . 
\label{c-bound}
\end{equation}
Initially $Z_0 = \{ \mathbf{0} \}$, where $\mathbf{0}$ denotes the origin.  
The ``monotonicity" condition we impose is that these rates are increasing in $S$:
\begin{equation}
\mbox{ 
if $v, v^\prime$ are adjacent to $S$ then 
$q(S, S \cup \{v\}) \le q(S\cup \{v^\prime\}, S \cup \{v, v^\prime\}) $
}.
\label{SZ-mono}
\end{equation}
Note that we do not assume any kind of spatial homogeneity.
\begin{Proposition}
Let $A$ be an arbitrary subset of vertices $\Ints^2 \setminus \{ \mathbf{0} \}$, and consider
$T:= \inf \{t: Z_t \cap A \mbox{ is non-empty.} \}$.
Under assumptions (\ref{c-bound}, \ref{SZ-mono}),
\[ \var \ T \le \Ex T/c_* . \]
\end{Proposition}
{\bf Proof.} 
Condition (\ref{SZ-mono}) allows us to couple versions 
$(Z^\prime_t,Z^{\prime \prime}_t)$ of the process starting from states 
$S^\prime \subset S^{\prime \prime}$, 
such that in the coupled process we have $Z^\prime_t \subseteq Z^{\prime \prime}_t$ for all $t \ge 0$.
In particular, 
$h(S):= \Ex_S T$ satisfies condition (\ref{def-monotone}). 
To deduce the result from Lemma \ref {L-simple-1} we need to show that, for any given possible transition 
$S_0 \to S_0 \cup \{v_0\}$, we have
\begin{equation}
h(S_0) \le h(S_0 \cup \{v_0\}) + 1/c_* . 
\label{hh7}
\end{equation}
Now by running the process started at $S_0$ until the first time $T^*$ this process contains $v_0$, 
and then coupling the future of that process to the process started at $S_0 \cup \{v_0\}$, we have
$h(S_0) \le \Ex _{S_0}T^* + h(S_0 \cup \{v_0\})$. 
And $\Ex _{S_0}T^* \le 1/c_*$ by (\ref{c-bound}), establishing (\ref{hh7}).

\subsection{A multigraph process}
\label{sec:multi}
The process here arises in a broad ``imperfectly observed networks" program described in  \cite{aldous-lisha}.
We give two examples  of the application of Lemma \ref{L-simple-1}  to this process.

Take a  finite connected graph $(\bV,\bE)$ with edge-weights $\bw = (w_e)$, where $w_e > 0 \ \forall e \in \bE$.
Define a multigraph-valued process as follows.  
Initially we have the vertex-set $\bV$ and no edges.
For each vertex-pair  $e = (vy) \in \bE$, edges $vy$ appear at the times of a Poisson 
(rate $w_e$) process, independent over $e \in \bE$.  
So at time $t$ the state of the process, $Z_t$ say, is a multigraph with $N_e(t) \ge 0$ copies of edge $e$, where 
$(N_e(t), \ e \in \bE)$ are independent Poisson($t w_e $) random variables.

We study how long until $Z_t$ has various connectivity properties.
Specifically, consider
\begin{itemize}
\item $\Tspan_k = \inf \{ t: Z_t \mbox{ contains $k$ edge-disjoint spanning trees} \} .$
\item $ \Ttria_k = \inf \{t: Z_t \mbox{ contains $k$ edge-disjoint triangles} \} $.
\end{itemize}
Here we regard the $N_e(t)$ copies of $e$ as disjoint edges.
Remarkably, Lemma \ref{L-simple-1} enables us to give  simple proofs of 
``weak concentration" bounds which do not depend on the underlying weighted graph.
\begin{Proposition}
\label{P:spantree}
\begin{eqnarray}
 \frac{ \mathrm{s.d.} (\Tspan_k)}{\Ex \Tspan_k} &\le& \frac{1}{\sqrt{k}}, \ k \ge 1 .  \label{bd-span}\\
\frac{ \mathrm{s.d.} (\Ttria_k)}{\Ex \Ttria_k} & \le  & \left(\frac{e}{e-1}\right)^{1/2} k^{-1/6} , \ k \ge 1 .  \label{bd-tria}
\end{eqnarray}
\end{Proposition}
Using the continuization device at (\ref{continuization}), the same bounds hold in the discrete-time model where edges $e$ arrive IID with 
probabilities proportional to $w_e$

We conjecture that some similar result holds for 
\[ T_k^\prime := \inf \{ t: Z_t \mbox{ is $k$-edge-connected} \} . \]
But proving this by our methods would require some structure theory 
(beyond Menger's theorem) for $k$-edge-connected graphs, and it is not clear whether 
relevant theory is known.

\paragraph{Proof of (\ref{bd-span}).} 
We will apply Lemma \ref{L-simple-1}.
Here the states $S$ are multigraphs over $\bV$, and $h(S)$ is the expectation, starting at $S$, 
of the time until the process contains $k$ edge-disjoint spanning trees.
What are the possible values of 
$h(S) - h(S \cup \{e\})$, 
where $S \cup \{e\}$ denotes the result of adding an extra copy of $e$ to the multigraph $S$?

Consider the ``min-cut" over proper subsets $S \subset \bV$
\[ \gamma := \min_S w(S,S^c)  \] 
where $w(S,S^c) = \sum_{v \in S, y \in S^c} w_{vy}$.
Because a spanning tree must have at least one edge across the min-cut,
\begin{equation}
\Ex \Tspan_k \ge k/\gamma. 
\label{gamma-1}
\end{equation}
On the other hand we claim
\begin{equation}
h(S) - h(S \cup \{e\}) \le 1/\gamma .
\label{gamma-2}
\end{equation}
To prove this, take the natural coupling $(Z_t,Z^+_t)$ of the processes started from  $S$ and from $S \cup \{e\}$, and run the coupled process until $Z^+_t$ contains $k$ edge-disjoint spanning trees.  
At this time, the process $Z_t$ either contains $k$ edge-disjoint spanning trees, or else 
contains $k-1$ spanning trees plus two trees 
(regard as edge-sets $\bt_1$ and $\bt_2$) such that $\bt_1 \cup \bt_2 \cup \{e\}$ is a spanning tree.  
So the extra time we need to run $(Z_t)$ 
is at most the time until some arriving edge links $\bt_1$ and $\bt_2$, which has mean at most $1/\gamma$.
This establishes (\ref{gamma-2}), and then  Lemma \ref{L-simple-1} establishes (\ref{bd-span}).

\medskip
\paragraph{Proof of (\ref{bd-tria}).} 
The key ingredient is a coupling construction which will establish
\begin{Lemma}
\label{L:triangles}
\[ \frac{\Ex \Ttria_1}{\Ex \Ttria_k} \le a(k) : = \inf_{0 < q < 1} \frac{q}{1 - (1-q^3)^k} . \]
\end{Lemma}
Granted Lemma \ref{L:triangles}, we start by repeating the format of the previous proof. 
To bound $h(S^\prime) - h(S)$ for a possible transition (one added edge) $S \to S^\prime$, 
we take the natural coupling of the processes started from $S$ and $S^\prime$, run until the latter 
contains $k$ edge-disjoint triangles.  
At that time the former process contains $k-1$ edge-disjoint triangles, so the mean extra time until it contains 
$k$ edge-disjoint triangles is bounded (crudely) by $\Ex \Ttria_1$.
Then Lemma \ref{L-simple-1} says
\[ \frac{\var  \ \Ttria_k}{\Ex  \Ttria_k} \le \Ex \Ttria_1 . \]
Combining with Lemma \ref{L:triangles},
\[ \frac{ \mathrm{s.d.} (\Ttria_k)}{\Ex \Ttria_k} \le \sqrt{a(k)} . \]
Because $1 - x \le e^{-x}$ we have
$\frac{q}{1 - (1-q^3)^k}  \le \frac{q}{1 - \exp(-q^3k)}$,
 so taking $q = k^{-1/3}$ gives the stated bound.
 
 \paragraph{Proof of Lemma \ref{L:triangles}.}
 Recall $(N_e(t), 0 \le t)$ is the rate-$w_e$ Poisson counting process for occurrence of copies of edge $e$.
 Consider the following coupling construction, done independently for different $e$. 
 Fix $0<q<1$.
 Delete edge copies with probability $1-q$, and write $(\widetilde{N}_e(t), 0 \le t)$ for the process of retained edges; then rescale time 
 by setting $N^*_e(t) = \widetilde{N}_e(t/q)$.
 So $(N_e(\cdot),N^*_e(\cdot))$ is a coupling of two dependent rate-$w_e$ Poisson processes.
 It is enough to show, writing $T$ for $\Ttria$,
 \begin{equation}
 \frac{\Ex T^*_1}{\Ex T_k} \le \frac{q}{1 - (1-q^3)^k} . 
 \label{TTkq}
 \end{equation}
 where $T^*_1$ refers to the $(N^*_e(\cdot), e \in \EE)$ process.  
 Now $T^*_1/q$ is the time of appearance of the first triangle in the  $(\widetilde{N}_e(\cdot), e \in \EE)$ process,
 and so $T^*_1/q \le \widetilde{T}$ where the latter is defined by
 \begin{quote}
 $ \widetilde{T} = T_k$ on the event $A = $ ``at least one of the $k$ triangles seen at $T_k$ in the $N^*(\cdot)$ process is retained in the 
 $ \widetilde{N}(\cdot)$ process";
 
 $ \widetilde{T} = T_k + Y$ otherwise, where $Y$ is the subsequent time until a triangle appears in the  $ \widetilde{N}(\cdot)$ process,
 if all edges are deleted at time $T_k$.
 \end{quote}
 But $\Ex Y = \Ex T_1/q  = \Ex T^*_1/q$, so
 \begin{eqnarray*}
 \Ex T^*_1/q   \le   \Ex \widetilde{T} &=&\Ex T_k + \Pr(A^c)\Ex Y \\
 &\le& \Ex T_k + (1-q^3)^k   \Ex T^*_1/q
 \end{eqnarray*}
 which rearranges to inequality (\ref{TTkq}).

 \subsection{Coverage processes}
 \label{sec:coverage}
 The topic of {\em coverage processes} is centered upon spatial or combinatorial variants of the coupon collector's problem; 
 see the monograph \cite{hall} and scattered examples in \cite{PCH}. 
 Classical theory concerns low-parameter models for which the cover time $T_n$ of a ``size $n$" model can be shown to have a limit
 distribution after scaling:
 $(T_n - a_n)/b_n \to_d \xi$ for explicit $a_n ,b_n$.
 In many settings, Lemma \ref{L-simple-1} can be used to give a weak concentration result for models with much less regular structure. 
 Here is a very simple example, whose one-line proof is left to the reader.
 \begin{Proposition}
 Let $G$ be an arbitrary $n$-vertex graph, and let $(V_i, i \ge 1)$ be IID uniform random vertices.  
 Let $T$ be the smallest $t$ such that every vertex is contained in, or adjacent to, the set 
 $\{V_i, 1 \le i \le t\}$.  Then 
 $\var \ T \le n \ \Ex T$.
 \end{Proposition}
 
 For a sequence $(G_n)$ of  sparse graphs, $\Ex T_n$ will be of order $n \log n$, so the bound says that 
 $\mathrm{s.d.}(T_n)/\Ex T = O(1/\sqrt{\log n})$.

\section{The FPP model}
\label{sec:introperc}
 As in section \ref{sec:multi}  we start with a finite connected graph $(\bV,\bE)$ with edge-weights $\bw = (w_e)$, where $w_e > 0 \ \forall e \in \bE$, but the model here is different.
 To the edges $e \in \bE$ attach independent Exponential(rate $w_e$) random variables $\xi_e$.
 For each pair of vertices $(v^\prime,v^{\prime\prime})$ there is a random variable $X(v^\prime,v^{\prime\prime})$ which can be 
 viewed in two equivalent ways:
 
 viewing $\xi_e$ as the length of edge $e$, then $X(v^\prime,v^{\prime\prime})$ is the length of the shortest route from $v$ to $v^\prime$;
 
  viewing $\xi_e$ as the time to traverse edge $e$, then $X(v^\prime,v^{\prime\prime})$ is the ``first passage percolation" time from $v$ to $v^\prime$.
  
  \medskip
  \noindent
  Taking the latter view, we call this the {\em FPP model}, and call $X(v^\prime,v^{\prime\prime})$ the {\em FPP time} 
  and $\xi_e$ the {\em traversal time}.
  This type of model and many generalizations have been studied extensively in several settings, in particular
  \begin{itemize}
  \item FPP with general IID weights on $\Ints^d$ \cite{kesten-FPP,50years}
  \item FPP on classical random graph ({E}rd{\H o}s-{R}\'enyi or configuration) models  \cite{shankar-1,shankar-2}
\item  and a much broader ``epidemics and rumors on complex networks"  literature \cite{draief,peng}.
  \end{itemize}
  However, this literature invariably starts by assuming some specific  graph model; we do not know any  ``general results" which relate properties of the FPP model to properties of a general underlying graph. 
 As an analogy, another structure that can be associated with a weighted finite graph is a finite reversible Markov chain; 
the established theory surrounding mixing times of Markov chains \cite{aldous-fill-2014,MCMT,monte} does contain ``general results" 
relating properties of the chain to properties of the  underlying graph. 

This article makes a modest start by studying the
 ``weak concentration" property: when it is true that $X(v^\prime,v^{\prime\prime})$ is close to its 
 expectation?
 We  can reformulate the FPP model as a set-valued process (see section \ref{sec:prelim}(a) for details) and then 
 Lemma \ref{L-simple-1} immediately implies
the following result.
  \begin{Proposition}
 \label{Pvar}
 $\var \ X(v^\prime,v^{\prime\prime}) \le \Ex X(v^\prime,v^{\prime\prime})/ w_*$ for 
 $w_* := \min \{w_e: \ e \in \bE\} $.
 \end{Proposition}
 Bounds of this type are classical on $\Ints^d$  \cite{kesten-speed}. 
 
Proposition \ref{Pvar} implies that on any unweighted graph ($w_e = 1$ for all edges $e$), 
the spread of $X = X(v^\prime,v^{\prime\prime})$ is at most order 
$\sqrt{ \Ex X}$.  
 For many specific graphs, stronger concentration results are known. 
 For $\Ints^2$ there is extensive literature (see \cite{chatterjee} for a recent overview) on the longstanding conjecture 
 that the spread is order $(\Ex X)^{1/3}$.
 For the complete graph (e.g. \cite{me-FMIE} sec. 7.3) and for sparse random graphs on $n$ vertices \cite{shankar-1} the spread of $X/\Ex X$ is typically 
 of order $1/\log n$.

In contrast, in this paper we study the completely general case 
where the edge-weights $w_e$ may vary widely over the different edges 
$e \in \bE$.
Here is our main result (conjectured in \cite{me-FMIE} sec. 7.4).
Given a pair $(v^\prime,v^{\prime\prime})$, there is a
random path $\bm{\pi}(v^\prime,v^{\prime\prime})$ that attains the FPP time $X(v^\prime,v^{\prime\prime})$.
Define $\Xi (v^\prime,v^{\prime\prime}) := \max \{ \xi_e: \ e \in \bm{\pi}(v^\prime,v^{\prime\prime}) \}$ 
as the maximum edge-traversal time in this minimal-time path.
 Recall the ``$L^0$ norm" 
 \[ || V ||_0 := \inf \{\delta: \Pr( |V| > \delta) \le \delta \} . \]
 \begin{Theorem}
 \label{Tmain}
 There exist functions $\psi^+$ and  $\psi_- : (0,1] \to (0,\infty)$ 
 such that $\psi^+(\delta) \downarrow 0$ as $\delta \downarrow 0$, and 
 $\psi_-(\delta) > 0$ for all $\delta > 0$, and
 such that, for all finite connected edge-weighted  graphs and all vertex pairs $(v^\prime,v^{\prime\prime})$,
 \[   \psi_- \left( \left| \left| \frac{\Xi(v^\prime,v^{\prime\prime}) }{\Ex X(v^\prime,v^{\prime\prime})} \right| \right|_0  \right)
 \le
  \frac{ \mathrm{s.d.}(X(v^\prime,v^{\prime\prime})) }{\Ex X(v^\prime,v^{\prime\prime})} 
 \le \psi^+ \left( \left| \left| \frac{\Xi(v^\prime,v^{\prime\prime}) }{\Ex X(v^\prime,v^{\prime\prime})} \right| \right|_0  \right) .
 \]
 \end{Theorem}
In words, $X/\Ex X$ has small spread if and only if $\Xi/\Ex X$ is small.
Intuition for this result comes from the ``almost disconnected" case where the path $\bm{\pi}(v^\prime,v^{\prime\prime})$ 
must contain a specific ``bridge" edge $e$ with small $w_e$; if $1/w_e$ is not $o(\Ex X)$ then the contribution to X from the 
traversal time $\xi_e$ is enough to show that $X$ cannot have weakly concentrated distribution.
 
 The proof of the upper bound starts from our ``more general" inequality, Lemma \ref{L-simple-2}.  
 See section \ref{sec:UB} for an outline.
 
Theorem \ref{Tmain} is unsatisfactory in that the conditions are not directly on the edge-weights $\bw = (w_e)$.
 By analogy with the bounds (\cite{monte} sec. 3.2) for Markov chain mixing times in terms of the {\em spectral profile},
it seems likely that Proposition \ref{Pvar} can be extended to give more widely applicable upper bounds on spread in terms of extremal flow rates 
$w(S,S^c) = \sum_{v \in S,y \in S^c} w_{vy}$.
However we do not have any conjecture for two-sided bounds analogous to Theorem \ref{Tmain}.

A particular case of FPP in our setting -- 
 edge-weights $w_e$  varying widely over the different edges -- 
 is studied in detail in recent work of  Chatterjee-Dey \cite{c-dey}. 
In their model $\bV = \Ints^d$ with $w_{xy} = ||y-x||^{- \alpha + o(1)}$.  
Results and conjectures from that paper are consistent with our Theorem \ref{Tmain}, 
which says that properties 
$X(\mathbf{0},nz)/\Ex X(\mathbf{0},nz) \to_p 1$
and 
$\Xi(\mathbf{0},nz)/\Ex X(\mathbf{0},nz) \to_p 1$
must either both hold or both fail.
They identify several qualitatively different regimes.  
In their {\em linear growth regime} ($\alpha > 2d+1$) 
they show $\Ex X(\mathbf{0},nz)/n$ converges to a nonzero constant 
and show that both properties hold.
 In their {\em super-linear growth regime} ($\alpha\in (2d,2d+1)$)
they show $X(\mathbf{0},nz) = n^{\alpha-2d+o(1)}$; here their analysis suggests both 
properties fail.
For $\alpha \in (d,2d)$ they show $X(\mathbf{0},nz)$ grows as a power of $\log n$, 
and their arguments suggest both properties hold.
The qualitative behavior and proof techniques in \cite{c-dey} are different in these different regimes, 
whereas our Theorem \ref{Tmain} is a single result covering all regimes, albeit a less explicit result.

 \subsection{Some preliminaries}
 \label{sec:prelim}
 \noindent
{\bf (a)}
We will view the FPP process started at $v^\prime$ as a process $(Z_t)$ taking values in the (finite) space of subsets $S \subseteq \bV$ of vertices, that is as
 \begin{equation}
  Z_t := \{v: X(v^\prime,v) \le t\} . 
  \label{def-Z}
  \end{equation}
 The assumption of Exponential distributions implies that $(Z_t)$ is the continuous-time Markov chain with 
 $Z_0 = \{v^\prime\}$ and transition rates
 \[ S \to S \cup \{y\}: \quad  \mbox{ rate }  w(S,y):= \sum_{s \in S} w_{sy} \quad (y \not\in S).
 \]
 So we are in the setting of Lemmas \ref{L-simple-1} and \ref{L-simple-2}.
 Given a target vertex $v^{\prime \prime}$ the FPP  time $X(v^\prime,v^{\prime \prime})$ is a stopping time 
 $T$ of the form (\ref{T-def}), and the function $h(S)$ at (\ref{h-def}) is just the function
  \[ h(S):= \Ex \min_{v \in S}  X(v,v^{\prime\prime}) .\]
  For a possible transition $S \to S \cup \{y\} $, 
 by considering the first time $y$ is reached in the process started from $S$ we have 
 $h(S) \le  \sfrac{1}{w(S,y)} + h(S \cup \{y\}) $ and so
 \[ h(S) -  h(S \cup \{y\}) \le  \sfrac{1}{w(S,y)} \le \sfrac{1}{w_*} \]
 for $w_* := \min \{w_e: \ e \in \bE\} $.  Now Lemma \ref{L-simple-1}  does indeed imply Proposition \ref{Pvar}, 
 as stated earlier.

\medskip \noindent
 {\bf (b)} FPP times such as $X(v^\prime,v^{\prime\prime})$ are examples of distributions $Y>0$ with the 
 {\em submultiplicative} property
 \[ \Pr (Y > y_1 + y_2) \le \Pr(Y > y_1) \ \Pr(Y > y_2); \quad y_1, y_2 > 0 .\]
We will need the general bound provided by the following straightforward lemma.
\begin{Lemma}
\label{L:subm}
There exists a function $\gamma(u) > 0$ with $\gamma(u) \downarrow 0$ as $u \downarrow 0$ such that,
for every  submultiplicative $Y$ and every event $A$,
\[ \frac{\Ex Y \ind_{A}}{\Ex Y} \le \gamma( \Pr(A)) . \]
\end{Lemma}
Note it follows that
\begin{equation}
 \Pr \left( Y \ge y\right) \ge \gamma^{-1} \left( \sfrac{\Ex Y \ind_{\{Y\ge y\}}}{\Ex Y} \right)  
 \label{gamma-inv}
 \end{equation}
for the inverse function
$\gamma^{-1}(u) \downarrow 0$ as $u \downarrow 0$.

\medskip \noindent
 {\bf (c)}
We will need an elementary stochastic calculus lemma.
\begin{Lemma}
\label{Lsc}
Let $T_1$ be a stopping time with (random) intensity $\eta_t$ -- that is, 
$\Pr(t < T_1< t +dt \vert \FF_t) = \eta_t \ dt$ on $\{T_1 > t\}$. 
Let $\zeta$ be another stopping time such that
$\eta_t \ge c$ on $\{\zeta > t\}$, for constant $c > 0$.
Then
\[ \Pr (T_1 \le \zeta \wedge t_0) \ge (1 - e^{-ct_0}) \Pr(\zeta > t_0) . \]
\end{Lemma}
{\bf Proof.} 
Applying the optional sampling theorem to the martingale 
$\exp ( \int_0^t \eta_s \ ds) \ \cdot \ind_{\{T>t\}}$ 
and the stopping time $\zeta \wedge t_0$ shows 
\[  \mbox{ $ 
\Ex [ \exp ( \int_0^{\zeta \wedge t_0} \eta_s \ ds) \ \cdot \ind_{\{T>\zeta \wedge t_0\}} ] = 1. 
$} \]
Now
\[ \mbox{ $ 
\exp ( \int_0^{\zeta \wedge t_0} \eta_s \ ds) 
\ge \exp ( c \ ( \zeta \wedge t_0)) \ge 1 + (e^{ct_0}-1) \ind_{\{\zeta > t_0\}} 
$}
\]
and so
\[ \Pr (T_1 > \zeta \wedge t_0) +  (e^{ct_0}-1) \Pr(T_1 > \zeta \wedge t_0, \zeta > t_0) \le 1 . \]
That is,
\begin{eqnarray*}
\Pr (T_1 \le \zeta \wedge t_0) &\ge& (e^{ct_0}-1) \Pr(T_1 > \zeta \wedge t_0, \zeta > t_0) \\
&\ge & (e^{ct_0}-1)  (\Pr(\zeta > t_0) - \Pr (T_1 \le \zeta \wedge t_0) )
\end{eqnarray*}
which rearranges to the stated inequality.

\section{Proof of upper bound in Theorem \ref{Tmain}}
\label{sec:UB}

 Fix $v^\prime$ and  as above view the FPP process started at $v^\prime$ as the continuous-time Markov chain $(Z_t)$ 
 at (\ref{def-Z}) on the space of subsets of $\bV$.
 Fix a target $v^{\prime\prime}  \ne v^\prime$ and in the following write $S$ for an arbitrary subset of vertices containing $v^\prime$.  
 Write
 \[ h(S):= \Ex \min_{v \in S}  X(v,v^{\prime\prime}) \]
 for the mean  percolation time from $S$ to $v^{\prime\prime}$, so $h(S) = 0$ iff $v^{\prime\prime} \in S$.
 So 
 \[ T := X(v^\prime,v^{\prime\prime})  = \inf \{t: h(Z_t) = 0\} \]
 is a stopping time for $(Z_t)$. 
 Note that $t \to h(Z_t)$ is decreasing and 
 $h(S) \le \Ex T$.
 
\medskip \noindent
{\bf Outline of proof.}
Step 1 is to translate Lemma \ref{L-simple-2} into our FPP setting; this shows it is enough to prove that in all transitions $S \to S \cup \{y\}$ the decrements
$h(S) - h(S \cup \{y\}) $ are  $o(\Ex T)$; so suppose not, that is suppose some are $\Omega (\Ex T)$.
Step 2 shows that  in some such transitions, the used edge $vy$ will (with non-vanishing probability) have traversal time $\xi_{vy}$ also of order $\Omega (\Ex T)$.
Step 3 shows that some such edges will (with non-vanishing probability) be in the minimal path.

\medskip \noindent
{\bf Details of proof.}
{\bf Step 1.}
Substituting $2 \delta$ for $\delta$, Lemma \ref{L-simple-2} says 
that for arbitrary $\delta, \eps > 0$
\begin{equation}
 \frac{\var \ T}{(\Ex T)^2} \le 2 \delta  + \eps
+\frac{  \Ex \int_0^T \ind_{\{q_\delta(Z_u) \ge \eps\}} du }{\Ex T} 
\label{vET}
\end{equation}
where
\[ q_\delta(S) := \sum_{y: \ h(S) - h(S \cup \{y\}) > 2\delta \Ex T} 
\ w(S,y) (h(S) - h(S \cup \{y\}) ) .\]

Informally, as  stated in the {\em outline of proof}, 
this shows it will suffice to prove that in transitions $S \to S \cup \{y\}$ the decrements
$h(S) - h(S \cup \{y\}) $ are  $o(\Ex T)$.

\bigskip \noindent
{\bf Step 2.}
 Now fix $0 < u_1 < u_2$.
\begin{Lemma}
\label{Lstep2}
Condition on $Z_{t_0} = S_0$.
The event
\begin{quote}
during $[t_0,t_0+u_2]$ the process $(Z_{t_0+u})$ makes a transition 
$S \to S \cup \{y\}$ such that $h(S) - h(S \cup \{y\}) > \delta \Ex T$ 
and using an edge $vy$ for which $\xi_{vy} > u_1$
\end{quote}
has probability at least 
$(1 - \sfrac{u_2}{\delta \Ex T})^+  \  \left( 1 - \exp(\sfrac{- (u_2-u_1)q(S_0)}{2 \Ex T} ) \right)$.
\end{Lemma}
{\bf Proof.} 
Define 
\[ \zeta = \inf\{t > t_0: h(Z_t) \le h(S_0) - \delta \Ex T\} \]
so that, because $t \to h(Z_t)$ is decreasing,
\begin{equation}
 h(S_0) - h(Z_{t_0 +u}) \le \delta \Ex T  \mbox{ on } \{\zeta > t_0 +u\} .
 \label{hShZ}
 \end{equation}
Define
\[ q^*(S) := \sum_{y: \ h(S) - h(S \cup \{y\}) > \delta \Ex T} 
\ w(S_0,y) (h(S) - h(S \cup \{y\}) ) .\]
Note this is the definition of $q_\delta(S)$ modified by changing the constraint ``$> 2 \delta \Ex T$" to 
``$>  \delta \Ex T$" and changing $w(S,y)$ to $w(S_0,y)$.  
Observe the relation, for $Z_{t_0+u} \supseteq S_0$,
{\small 
\begin{eqnarray}
h(Z_{t_0+u}) - h(Z_{t_0+u} \cup \{y\}) \! \!& \ge &
h(Z_{t_0+u}) - h(S_0 \cup \{y\}) \nonumber \\
&=& \! \! \!  [ h(S_0) - h(S_0 \cup \{y\}) ] \! - \! [ h(S_0) - h(Z_{t_0+u})] .
\label{hhhh}
\end{eqnarray}
}
For $y$ satisfying the constraint
\[ h(S_0) - h(S_0 \cup \{y\}) > 2 \delta \Ex T \]
in the definition of $q(S_0)$, 
and on the event $\{\zeta > t_0 +u\}$,
inequalities (\ref{hShZ},\ref{hhhh})  show that $y$ satisfies the constraint in the definition of $q^*(Z_{t_0+u})$, 
and also show that 
\[ h(Z_{t_0+u}) - h(Z_{t_0+u} \cup \{y\}) \ge \sfrac{1}{2} [ h(S_0) - h(S_0 \cup \{y\}) ] , \]
implying that 
\begin{equation}
 q^*(Z_{t_0+u}) \ge \sfrac{1}{2} q(S_0) \mbox{ on } \{\zeta > t_0 +u\} .
\label{qZtu}
\end{equation}
For $S \supseteq S_0$ we have, from the definition of $q^*(S)$,  a crude bound
\[ q^*(S) \le  (\Ex T ) 
 \sum_{y: \ h(S) - h(S \cup \{y\}) > \delta \Ex T} 
\ w(S_0,y)  \]
and applying this to the left side of (\ref{qZtu}) gives
\[
 \sum_{y: \ h(Z_{t_0+u}) - h(Z_{t_0+u} \cup \{y\}) > \delta \Ex T} 
\ w(S_0,y) \ge \frac{q(S_0)}{2 \Ex T}  \mbox{ on } \{\zeta > t_0 +u\} .
\]
Over the time interval $[t_0 + u_1,t_0 + u_2]$ the left side is the intensity of an event which implies the event
($D$, say) in Lemma \ref{Lstep2}   
(in particular, an edge $vy$ is used with $v \in S_0 = Z(t_0)$, and so its ``age" $\xi_{vy}$ must be at least $u_1$).
We are now in the setting of Lemma \ref{Lsc}, which  shows that
\[ 
\Pr(D) \ge \Pr (\zeta > t_0 +u) \ \left( 1 - \exp(\sfrac{- (u_2-u_1)q(S_0)}{2 \Ex T} ) \right) 
. \]
By the martingale property (\ref{MgT}) of $h(Z_t) + t$ and Markov's inequality 
\begin{equation}
 \Pr (\zeta \le t_0 +u)  = \Pr( h(Z_{t_0 +u_2}) \le h(S_0) - \delta \Ex T) \le \sfrac{u_2}{\delta \Ex T}
\label{hZS}
\end{equation}
establishing the bound stated
in Lemma \ref{Lstep2}.

\bigskip \noindent
 {\bf Step 3.}
 \begin{Lemma}
 \label{Lstep3}
 Conditional on the process $(Z_u)$ making at time $t_0$ a transition  
 $S \to S \cup \{y\}$ using an edge $vy$, the probability that edge $vy$ is in the minimal 
 path $\bm{\pi}(v^\prime,v^{\prime\prime}) $ is at least 
 $\gamma^{-1}(\sfrac{h(S) - h(S \cup \{y\})}{\Ex T})$, for the inverse function 
 $\gamma^{-1}(\cdot) > 0$ at (\ref{gamma-inv}).
 \end{Lemma}
 {\bf Proof.}
 Condition as stated.  So after the transition we have
 \[ T - t_0 = \min(A, B) \]
 where $B = X(y,v^{\prime \prime})$ and $A = \min_{v \in S} X^*(v,v^{\prime \prime})$ where $X^*$ denotes the 
 minimum over paths not using $y$.  
 The probability in question equals $\Pr(B < A)$; note that $A$ and $B$ are typically dependent.
 
 It is easy to check that the solution to the problem
 \begin{quote}
 given the distribution of $A > 0$, construct a r.v. $B^\prime >0$ to minimize $\Pr(B^\prime <A)$ subject to the constraint that 
 $\Ex A - \Ex \min(A,B^\prime )$ takes a given value
 \end{quote}
 is to take $B^\prime = A \ind_{\{A \le a_0\}}$ 
 for $a_0$ chosen to satisfy the constraint.
 So in our setting,
 \[ \Pr(B < A) \ge \Pr(A > a_0) \]
for $a_0$ chosen to satisfy 
\[  \Ex A \ind_{\{A > a_0\}} = \Ex A - \Ex \min(A,B ) = h(S) - h(S \cup \{y\}) .\]
Because $A$ has the submultiplicative property,
Lemma \ref{L:subm} shows that 
\[ \Pr(A > a_0) \ge \gamma^{-1} \left( \sfrac{h(S) - h(S \cup \{y\})  }{\Ex A} \right).
 \]
Because $\gamma^{-1}(\cdot)$ is decreasing and $\Ex A \le \Ex T$, we have established the lemma.

\bigskip \noindent
{\bf Step 4.}  We now combine the  ingredients above.
Define
$U_1 = \inf \{t \geq 0: q(Z_t) \ge \eps\}$ and inductively
\[ U_{j+1} = \inf \{t \geq U_j + u_2: q(Z_t) \ge \eps\}. \]
Note that the condition
\begin{equation}
j u_2  \le \int_0^T \ind_{\{q(Z_u) \ge \eps\}} du 
\label{ju2}
\end{equation}
is sufficient to ensure $U_j < \infty$. 
Condition on $U_j = t_0, Z_{U_j} = S_0$ and apply Lemmas \ref{Lstep2} and \ref{Lstep3}.
We deduce that, with probability at least
\[ (1 - \sfrac{u_2}{\delta \Ex T}) \  \left( 1 - \exp(\sfrac{- (u_2-u_1) \eps }{2 \Ex T} ) \right) 
\ \gamma^{-1}(\delta) \]
the minimal path  $\bm{\pi}(v^\prime,v^{\prime\prime}) $ contains an edge $vy$ with $\xi_{vy} \ge u_1$ and such that $y$ is first reached 
during $(U_j, U_{j+1}]$.
The latter property ensures these edges are distinct as $j$ varies (but note the corresponding $\xi_{vy}$ are dependent).
Summing over $j$, applying (\ref{ju2})  and taking expectation, 
we find that 
\[ N(u_1) := |\{e \in   \bm{\pi}(v^\prime,v^{\prime\prime}) : \ \xi_e \ge u_1\}| \]
satisfies
\begin{equation}
\Ex N(u_1) \ge (1 - \sfrac{u_2}{\delta \Ex T})^+ \  \left( 1 - \exp(\sfrac{- (u_2-u_1) \eps }{2 \Ex T} ) \right) 
\ \gamma^{-1}(\delta) \ \sfrac{1}{u_2} \Ex \int_0^T \ind_{\{q(Z_u) \ge \eps\}} du .
\label{Nuu}
\end{equation}
Writing $\Xi = \Xi(v^\prime,v^{\prime\prime})$, the event 
$\{N(u_1) \ge 1\}$ is the event $\{\Xi \ge u_1\}$. 
Also $u_1 N(u_1) \le T$, so we can write
$N(u_1) \le \frac{T}{u_1} \ind_{\{\Xi \ge u_1\}}$
and then
\[ \Ex N(u_1) \le \sfrac{\Ex T}{u_1} \ 
 \sfrac{ \Ex T \ind_{\{\Xi \ge u_1\}} }{\Ex T} 
  \le \sfrac{\Ex T}{u_1} 
\gamma (\Pr( \Xi > u_1) ) , \]
the final inequality by Lemma \ref {L:subm}.
Combining with (\ref{Nuu}) and rearranging gives a bound for the final term of (\ref{vET}):
\[ \frac{  \Ex \int_0^T \ind_{\{q(Z_u) \ge \eps\}} du }{\Ex T} 
\le 
\frac{u_2  \ \gamma (\Pr( \Xi > u_1) )    }  { u_1  (1 - \sfrac{u_2}{\delta \Ex T})^+ \  \left( 1 - \exp(\sfrac{- (u_2-u_1) \eps }{2 \Ex T} ) \right) 
\ \gamma^{-1}(\delta)  }  .
\]
This becomes a little easier to understand when we set $u_1 = u \Ex T$ and $u_2 = 2 u_1$; then (\ref{vET}) becomes
\begin{equation}
  \frac{\var \  T}{(\Ex T)^2} \le 2 \delta  + \eps + 
\frac{2  \ \gamma (\Pr( \sfrac{\Xi }{\Ex T}> u) )    }  {  (1 - \sfrac{u}{\delta })^+ \  \left( 1 - \exp(\sfrac{- u \eps }{2 } ) \right) 
\ \gamma^{-1}(\delta)  }  .
\label{fff}
\end{equation}
We could manipulate this into a complicated expression for 
an upper bound  function $\psi^+(\cdot)$ in Theorem \ref{Tmain}.
But it is simpler to observe that existence of {\em some} upper bound function is equivalent to 
the following limit assertion: 
for weighted graphs $\bw^{(n)}$ and FPP times 
$T^{(n)} = X^{(n)}(v^{(n)\prime},v^{(n)\prime\prime})$
\[ \mbox{ if } 
\frac{ \Xi^{(n)}(v^{(n)\prime},v^{(n)\prime\prime}) }{   \Ex T^{(n)}} 
\to_p 0 
\mbox{ then }  \frac{\var \  T^{(n)} }{(\Ex T^{(n)} )^2} \to 0 .
\]
And this is immediate from (\ref{fff}).

 \section{Proof of lower bound in Theorem \ref{Tmain}}
 \label{sec:LB}
 We described this lower bound assertion in \cite{me-FMIE} sec. 7.4 as ``intuitively clear (and not hard to prove)".
 The {\em intuition} is that, given a typical realization of the process where $X = X(v^\prime,v^{\prime\prime})$ has some value $X_0$, and some value $\Xi_0$ attained by some $\xi_e$, 
 there are other realizations for which this $\xi_e$ ranges  over the interval $[0,\Xi_0]$, for which therefore 
 the value of $X$ ranges over $[X_0 - \Xi_0, X_0]$, so the order of  s.d.($X$) should be at least 
 the order of $\Xi$.
 Formalizing this intuition is not quite trivial. 
 The proof below is rather crude -- likely there is a more elegant
 argument giving a better bound.
 
 \medskip \noindent
 {\bf Proof.}
 Fix a pair $(v^\prime,v^{\prime\prime})$ and write $X=    X(v^\prime,v^{\prime\prime})$
 and  $\Xi =   \Xi(v^\prime,v^{\prime\prime})$. 
 So $X$ is a function of the traversal times 
 $(\xi_e)$.
 If we couple $(\xi_e, e \in \bE)$ with a copy $(\xi^\prime_e, e \in \bE)$  of  $(\xi_e)$, then
 by an elementary inequality
 \begin{equation}
  \var \ X \ge \sfrac{1}{4} \Ex (X^\prime - X)^2 . 
  \label{var-eq}
  \end{equation}
 We will use the following coupling.
 Fix $0 < a<b$ and write $\EE(\lambda;a,b)$ for the distribution of an Exponential($\lambda$) random variable conditioned to lie in $[a,b]$.
 Define $(\xi^\prime_e, \xi_e)$ to be independent as $e$ varies, with
 
 if $\xi_e \notin [a,b]$ then $\xi^\prime_e  = \xi_e$
 
  if $\xi_e \in [a,b]$ then $\xi^\prime_e$ is conditionally independent of $\xi_e$ with conditional 
  distribution $\EE(w_e;a,b)$.
  
  \noindent
  Consider the set  $D_{ab}$ of edges $e$ in the minimal path $\bm{\pi} = \bm{\pi}(v^\prime,v^{\prime\prime})$ for which 
  $\xi_e \in [a,b]$.
  By considering the same path in the coupled copy,
  \begin{equation}
   X^\prime - X \le \sum_{e \in \bm{\pi}} (\xi^\prime_e - \xi_e)
   = \sum_{e \in D_{ab}} (\xi^\prime_e - \xi_e) .
    \label{XXe}
    \end{equation}
 To analyze this expression we will use the following lemma. 
 For $k \ge 1$ and $s > 0$ define
 \begin{equation}
  F_k(s) := \Ex  \left(  \max(0, s - \sum_{i=1}^k U_i) \ \right)^2 
  \label{eq-Fk}
  \end{equation}
 where the $U_i$ are independent Uniform$(0,1)$.
 \begin{Lemma}
 Let $(V_i, 1 \le i \le k)$ be independent $\EE(\lambda_i;a,b)$, for arbitrary $(\lambda_i)$. 
 Then for arbitrary $v_i \in [a,b]$,
 \begin{equation}
 \Ex  \left(  \max\left(0, \sum_{i=1}^k v_i - \sum_{i=1}^k V_i\right) \ \right)^2 
 \ge (b-a)^2 F_k\left(\sum_{i=1}^k \sfrac{v_i - a}{b-a} \right) .
 \label{baF}
 \end{equation}
 \end{Lemma}
 {\bf Proof.}  The distributions $\EE(\lambda_i;a,b)$ are stochastically decreasing with $\lambda$, so the 
 left side is minimized in the $\lambda_i \to 0$ limit, which is the uniform distribution on $[a,b]$; the limit value is obtained from (\ref{eq-Fk}) by scaling.

\medskip \noindent
We can now combine  inequalities (\ref{XXe}) and (\ref{baF}) to deduce that,
on the event $\{ |D_{ab}| = k\}$ for $k  \ge 1$,
  \[
 \Ex ( (X - X^\prime)^2 \vert \xi_e, e \in \bE) 
 \ge (b-a)^2 F_k\left(\sum_{e \in D_{ab}} \sfrac{\xi_e - a}{b-a} \right) .
 \]
 So
 \begin{equation}
   \Ex (X - X^\prime)^2 \ge (b-a)^2 \sum_{k \ge 1} 
 \Ex \left[ \ind_{\{|D_{ab}| = k\}} \  F_k\left(\sum_{e \in D_{ab}} \sfrac{\xi_e - a}{b-a} \right)\ \right] .
 \label{XXba}
 \end{equation}
 We need to lower bound the right side.  The issue is that, if $w_e$ is large then $\xi_e - a$ may be small.
 We handle this issue by considering two different values of $a$.
 Note that the functions $F_k(s)$ are decreasing in $k$ and increasing in $s$.
 Fix $a_1 < a_2 < b$, so that for $1 \le k \le K$
 \begin{eqnarray*}
  F_k\left(\sum_{e \in D_{a_1b}} \sfrac{\xi_e - a_1}{b-a_1} \right) & \ge & 
   F_K\left(\sum_{e \in D_{a_1b}} \sfrac{\xi_e - a_1}{b-a_1} \right)  \quad \mbox{(decreasing in $k$)} \\
  &\ge & 
 F_K\left(\sum_{e \in D_{a_2b}} \sfrac{\xi_e - a_1}{b-a_1} \right)  \quad \mbox{($D_{a_2b} \subseteq D_{a_1b}$)} \\
  & \ge & 
F_K\left( \ |D_{a_2b}| \ \sfrac{a_2 - a_1}{b-a_1} \right)  \quad \mbox{(increasing in $s$)} \\
& \ge &
F_K\left( \ \sfrac{a_2 - a_1}{b-a_1} \right)  \ \ind_{\{ |D_{a_2b}| \ge 1  \} } \quad \mbox{(increasing in $s$)} .
 \end{eqnarray*}
 Applying (\ref{XXba})  with $a = a_1$ and restricting the sum to $1 \le k \le K$,
  \begin{equation}
    \Ex (X - X^\prime)^2 \ge (b-a_1)^2 \ F_K\left( \ \sfrac{a_2 - a_1}{b-a_1} \right) 
  \ \Pr( |D_{a_1b}| \le K, \ |D_{a_2b}| \ge 1) .
  \label{XXba2}
  \end{equation}
 Because $a_1 |D_{a_1b}| \le X$, Markov's inequality tells us
$\Pr ( |D_{a_1b}| \ge K) \le \sfrac{1}{a_1K} \Ex X $.
Also the event $\{  |D_{a_2b}| \ge 1\}$ contains the event 
$\{a_1 \le \Xi \le b\}$, so 
\[ \Pr( |D_{a_1b}| \le K, \ |D_{a_2b}| \ge 1) 
\ge \Pr( \Xi \ge a_1) - \Pr( \Xi \ge b) - \sfrac{1}{a_1K} \Ex X . \]
Because $\Xi \le X$ we have 
$\Pr(\Xi \ge b) \le \Ex X/b$, and combining with (\ref{XXba2}) and (\ref{var-eq})
we conclude
\[ \var \ X \ge \sfrac{1}{4} (b-a_1)^2 \ F_K\left( \ \sfrac{a_2 - a_1}{b-a_1} \right) \ 
\left( \Pr( \Xi \ge a_1) - (\sfrac{1}{b} + \sfrac{1}{a_1K} )  \Ex X \right)^+ .
\]
This holds for arbitrary $0 < a_1 < a_2 < b$ and $K \ge 1$.
We now take $1 > \delta > 0$ and choose 
\[ a_1 = \delta \Ex X; \quad a_2 = 2a_1; \quad b = 3 \Ex X/\delta   \]
and find
\[ \frac{\var \ X}{(\Ex X)^2} \ge  \sfrac{1}{4} (\sfrac{3}{\delta} - \delta)^2 \ F_K\left( \ \sfrac{\delta^2}{3 - \delta^2} \right) 
\left( \Pr(\sfrac{ \Xi}{\Ex X} \ge \delta) - (\sfrac{\delta}{3} + \sfrac{1}{\delta K} )  \right)^+ .
\]
So finally choosing $K = K(\delta) \ge 3\delta^{-2}$, we see that the lower bound in Theorem \ref{Tmain} holds for
\[ \psi_-(\delta):= \sqrt{ 
 \sfrac{1}{4} (\sfrac{3}{\delta} - \delta)^2 \ F_K ( \ \sfrac{\delta^2}{3 - \delta^2} ) \ \sfrac{\delta}{3}
} . \]

 \section{Final remarks}
 \label{sec:remarks}
 {\bf (a).} The key point of results such as Proposition \ref{P:spantree} and Theorem \ref{Tmain} is that the bounds do 
 not depend on the size $n$ of the graph -- any 
 simple use of the method of bounded differences in these models would give bounds that did depend on $n$.

 \medskip \noindent
{\bf (b).}  A somewhat different general approach to  proving weak concentration for general coverage processes, 
assuming IID random subsets, was given in \cite{me-cover}, and used 
to prove  weak concentration for Markov chain cover times.
The latter result does not seem to follow easily from the methods in this paper.

\medskip \noindent
 {\bf (c).}  Theorem \ref{Tmain} went beyond Proposition \ref{Pvar} by using Lemma \ref{L-simple-2} instead of the simpler Lemma \ref{L-simple-1}, 
 and so one can imagine analogous improvements of the kind of results in sections \ref{sec:growth} - \ref{sec:coverage}.
 
\medskip \noindent
 {\bf (d).}  
 There is an analog of our FPP result for bond percolation.  
Let edges $e$ become ``open" at the independent random Exponential times $\xi_e$ of  rates $w_e$. 
Using Lemma \ref{L-simple-2}, it is proved in \cite{me-incipient} that 
(under minimal assumptions) the time at which a giant open component starts to emerge is weakly concentrated around its mean.

 \medskip \noindent
 {\bf (e).} 
 Presumably some analog of Theorem \ref{Tmain} remains true when the Exponential assumption is relaxed to 
 some much weaker ``shape of distribution" assumption, but we have not investigated that setting.

 \paragraph{Acknowledgements.}  I thank Jian Ding for extensive discussions and for dismissing earlier conjectures, 
 and Partha Dey for elucidating the connection with \cite{c-dey}.


 \end{document}